\documentclass[11pt,leqno]{article}

\usepackage{amsmath,amssymb,amsthm}
\usepackage{fullpage}
\usepackage{hyperref}
\usepackage{setspace}
\usepackage{dsfont}

\newtheorem{theorem}{Theorem}[section]
\newtheorem{lemma}[theorem]{Lemma}
\newtheorem{proposition}[theorem]{Proposition}

\newtheorem{remark}[theorem]{Remark}

\theoremstyle{remark}

\numberwithin{equation}{section}

\newcommand{\ls}{\leqslant}
\newcommand{\gr}{\geqslant}

\DeclareMathOperator{\Cov}{Cov}
\DeclareMathOperator{\Var}{Var}
\newcommand{\Hess}{\Sigma}
\hypersetup{colorlinks=true,linkcolor=blue,citecolor=blue,urlcolor=blue}

\begin{document}
\small

\title{\bf Wet-part doubling for log-concave probability\\ measures}
\author{Apostolos Giannopoulos and Minas Pafis}
\date{}
\maketitle

\begin{abstract}\footnotesize 
Let $\mu$ be a log-concave probability measure on $\mathbb R^n$, $n\gr 2$. We study the growth of $w_\mu(t)=\mu(\{x:q_{\mu}(x)\ls t\})$,
where $q_{\mu}$ is the Tukey half-space depth function of $\mu$. We prove that $w_{\mu}$ satisfies the doubling bound
$$w_\mu(\lambda t)\ls C^n\lambda\left(1+\frac{\log\lambda}{n-1}\right)^{n-1}w_\mu(t)$$
for every $\lambda>1$ and $0<t<1/\lambda$, where $C>0$ is an absolute constant. It follows that for every $\varepsilon>0$ we may write $w_\mu(\lambda t) \ls A_\varepsilon^n \lambda^{1+\varepsilon}w_\mu(t)$ for some constant $A_\varepsilon>0$. We construct a family of symmetric log-concave measures in every dimension $n\gr2$ showing that no linear doubling bound with a dimension-dependent constant holds uniformly over this class. The proof
of the doubling bound uses a pointwise comparison of half-space depth
with the Cram\'{e}r transform, dual radial coordinates, and a relative entropy estimate for transverse covariance determinants.
\end{abstract}

\section{Introduction}\label{sec:introduction}

The wet part of a convex body is the part of the body that can be removed by a half-space of small relative volume. It is a basic object in the theory of floating bodies, random polytopes, and approximation by polytopes. We begin with the classical results for convex bodies and then formulate the corresponding problem for log-concave probability measures.

Let $K$ be a convex body in $\mathbb{R}^n$ and let $\mu_K$ denote the normalized Lebesgue measure on $K$. Define
$$q_K(x)=\inf\{\mu_K(K\cap H): H\text{ is a closed half-space containing }x\},\quad x\in K,$$
and set $U_t(K)=\{x\in K:q_K(x)\ls t\}$ and $w_K(t)=\mu_K(U_t(K))$. The complement of the wet part $U_t(K)$ is the
interior of the convex floating body of $K$ at level $t$. Basic references on floating bodies are the articles of Sch\"{u}tt and Werner~\cite{Schutt-Werner-1990}, Nagy, Sch\"{u}tt and Werner~\cite{Nagy-Schutt-Werner-2019}, and the survey of Werner~\cite{Werner-2022}.

A fundamental result in this area is the economic cap-covering theorem of B\'{a}r\'{a}ny and Larman~\cite[Theorem~6]{Barany-Larman-1988}
which states that for sufficiently small $t$, there exist caps $C_1,\ldots,C_M$ of $K$ with $\mu_K(C_i)\ls A_nt$, and disjoint measurable 
subsets $C'_i\subseteq C_i$ of measure at least $a_nt$, such that
$$\bigcup_{i=1}^M C'_i\subseteq U_t(K)\subseteq\bigcup_{i=1}^M C_i.$$
In particular, $M\approx_n w_K(t)/t$. The proof of this theorem uses Macbeath's centrally symmetric regions \cite{Macbeath-1952} 
and develops the economical covering method of Ewald, Larman and Rogers~\cite{Ewald-Larman-Rogers-1970}. The classical theory also provides 
the containment of suitably small caps in enlarged covering regions. An account of the method is given in B\'{a}r\'{a}ny's survey article~\cite{Barany-2000}.

The starting point of our work is the wet-part comparison theorem of B\'{a}r\'{a}ny and Larman~\cite[Theorem~7]{Barany-Larman-1988}.
For $n\gr2$, for every $\lambda\gr1$, for sufficiently small $t>0$ and any convex body $K$ in $\mathbb{R}^n$, one has
\begin{equation}\label{eq:classical-doubling}
w_K(\lambda t)\ls C_n\lambda^n w_K(t),
\end{equation}
with a constant depending only on the dimension. This is a uniform power bound for the growth of wet parts, but the exponent $n$ is not
optimal. In this work we establish a substantially sharper bound for the wider class of log-concave probability measures.

The geometric significance of $w_K$ is especially clear from its applications to random polytopes. Let $P_N$ be the convex hull of
$N$ independent uniform points in $K$. The missing volume theorem of B\'{a}r\'{a}ny and Larman~\cite[Theorem~1]{Barany-Larman-1988}
states that
\begin{equation}\label{eq:intro-classical-random}
\mathbb E\,\mu_K(K\setminus P_N)\approx_n w_K(1/N)
\end{equation}
if $N$ is sufficiently large depending only on $n$. A classical identity of Efron~\cite{Efron-1965} asserts that 
$\mathbb E f_0(P_N)=N\,\mathbb E\mu_K(K\setminus P_{N-1})$, so the same geometric quantity also 
controls the expected number of vertices of $P_N$. Earlier work of R\'{e}nyi and Sulanke \cite{Renyi-Sulanke-1963}
concerned the planar case. B\'{a}r\'{a}ny~\cite{Barany-1989} subsequently studied intrinsic volumes and face numbers of $P_N$. 
Further work on random convex hulls and the combinatorics of their faces was done by B\'{a}r\'{a}ny and F\"uredi~\cite{Barany-Furedi-1988} and
Reitzner~\cite{Reitzner-2005}. See B\'{a}r\'{a}ny~\cite{Barany-2008} and Reitzner~\cite{Reitzner-Survey} for further results and references.

The behavior of $w_K(t)$ reflects the boundary geometry of $K$. The general convex body estimate is $w_K(t)\ls C_nt^{2/(n+1)}$ for
small $t$, with the exponent realized by balls and smooth convex bodies. On the other hand, if $P$ is an $n$-dimensional polytope, then
\begin{equation}\label{eq:intro-polytope-wet}
w_P(t)\approx_P t\bigl(\log(1/t)\bigr)^{n-1}
\end{equation}
as $t\to 0^+$. For the relation of $w_K$ with floating bodies and affine surface area we refer to \cite{Barany-Larman-1988, Schutt-Werner-1990, Werner-2022}. 

We consider the class of log-concave probability measures. Tukey's half-space depth is the natural
analogue of $q_K$. For a log-concave probability measure $\mu$, we define
$$q_\mu(x)=\inf\{\mu(H):H\text{ is a closed half-space containing }x\}.$$
The wet part of $\mu$ and its measure are
$$W_\mu(t)=\{x:q_\mu(x)\ls t\},\quad w_\mu(t)=\mu(W_\mu(t)).$$
Thus $w_\mu$ is the distribution function of half-space depth under $\mu$. For $\mu=\mu_K$ it agrees with $w_K$. 

Half-space depth was introduced by Tukey in~\cite{Tukey-1975}, and was further developed in robust statistics by Donoho and Gasko~\cite{Donoho-Gasko-1992}
and others~\cite{Rousseeuw-Ruts-1999,Liu-Parelius-Singh-1999,Zuo-Serfling-2000}. 

B\'{a}r\'{a}ny, Fradelizi, Goaoc, Hubard and Rote \cite{Barany-Et-Al-2020} proved that for arbitrary probability distributions 
the lower wet-part bound in~\eqref{eq:intro-classical-random} continues to hold but a logarithmic enlargement of the depth parameter 
is in general needed for the upper bound,  and they gave an example showing that this loss cannot be avoided. They also posed the 
question whether, for any $\lambda>1$, a uniform comparison between $w_\mu(t)$ and $w_\mu(\lambda t)$ holds true for log-concave probability measures.

We provide a strong affirmative answer to this question.

\begin{theorem}\label{th:near-linear}
There exists an absolute constant $C\gr1$ such that, for every log-concave probability measure $\mu$ on
$\mathbb R^n$, $n\gr2$, we hav
\begin{equation}\label{eq:near-linear-log}
w_\mu(\lambda t) \ls C^n\lambda\left(1+\frac{\log\lambda}{n-1}\right)^{n-1}w_\mu(t).
\end{equation}
for every $\lambda>1$ and every $0<t<1/\lambda$. Consequently, for every $0<\varepsilon\ls1$,
\begin{equation}\label{eq:near-linear-epsilon}
w_\mu(\lambda t) \ls C_1^n\varepsilon^{-(n-1)}\lambda^{1+\varepsilon}w_\mu(t).
\end{equation}
In particular, for any $\varepsilon>0$ there exists a constant $A_\varepsilon >0$, which is independent of $n$, such that
$w_\mu(\lambda t)\ls A_\varepsilon^n\lambda^{1+\varepsilon}w_\mu(t)$. 
\end{theorem}

The exponent $1$ cannot itself be used uniformly over all symmetric log-concave probability measures, even if we allow a constant depending on $n$.

\begin{theorem}\label{th:general-failure}
Let $n \gr 2$. For every $C_n>0$ there exists a symmetric log-concave probability measure $\mu$ on $\mathbb{R}^n$, and there exist
$\lambda>1$ and $0<t<1/\lambda$ such that 
$$w_\mu(\lambda t)>C_n \lambda w_\mu(t).$$
\end{theorem}

To prove Theorem~\ref{th:general-failure} we construct symmetric log-concave distributions $\mu_{n,a}$ and  pairs $(\lambda_{n,a},t_{n,a})$ for which
\begin{equation}\label{eq:intro-defect}
\frac{w_{\mu_{n,a}}(\lambda_{n,a}t_{n,a})}{\lambda_{n,a}w_{\mu_{n,a}}(t_{n,a})}\gr c_n\left(\frac{a}{\log a}\right)^{n-1}
\rightarrow\infty\quad(a\to\infty).
\end{equation}
The measures vary with $a$. This shows that we cannot have a uniform linear bound. However, the construction does not ensure that the 
logarithmic correction in Theorem~\ref{th:near-linear} is optimal.
Note that if $n=1$ then a simple calculation shows that $w_\mu(t)=\min\{2t,1\}$, so linear growth holds.

\smallskip

The proof of Theorem~\ref{th:near-linear} is based on the log-Laplace transform and its Legendre transform. We define
$$\Lambda(\theta)=\log\int_{\mathbb R^n}e^{\langle\theta,z\rangle}\,d\mu(z),\quad
\Lambda^*(x)=\sup_\theta\{\langle\theta,x\rangle-\Lambda(\theta)\},$$
where $\Omega=\{\theta:\Lambda(\theta)<\infty\}$, and for $\theta \in \Omega$ we set
$$x(\theta)=\nabla\Lambda(\theta),\quad \Hess(\theta)=\nabla^2\Lambda(\theta),\quad F(\theta)=\Lambda^*(x(\theta)),
\quad J(\theta)=\sqrt{\theta^T\Hess(\theta)\theta}.$$
The tilted probability measure
$d\mu_\theta(z)=e^{\langle\theta,z\rangle-\Lambda(\theta)}d\mu(z)$ has mean $x(\theta)$ and covariance $\Hess(\theta)$.
The next theorem is our starting point.

\begin{theorem}\label{th:pointwise}
There exist absolute constants $C_0 \gr c_0>0$ such that, for every log-concave probability measure $\mu$,
every $\theta\in\Omega$, and every $x=\nabla\Lambda(\theta)$,
\begin{equation}\label{eq:pointwise}
c_0\frac{e^{-F(\theta)}}{1+J(\theta)} \ls q_\mu(x)\ls C_0\frac{e^{-F(\theta)}}{1+J(\theta)}.
\end{equation}
\end{theorem}

For the upper bound in Theorem~\ref{th:pointwise} we evaluate the half-space with normal vector $\theta$. For the lower bound, we compare an
arbitrary supporting half-space through $x(\theta)$ with an appropriately tilted one. Relative entropy controls both the
change in the rate function $F$ and the one-dimensional variance factor.

Estimates for the relation of half-space depth of log-concave measures with thresholds for the expected measure
of random polytopes can be found in~\cite{BGP-2024-depth,BGP-2024-threshold}. Sharp estimates in the spirit of
Theorem~\ref{th:pointwise} can be found in the work of Brazitikos and Chasapis~\cite[Theorem~1.1]{BC-2025} who
showed that for every log-concave probability measure $\mu$ on $\mathbb{R}^n$, and every $\varepsilon \in (0,1)$, one has
\begin{equation}\label{eq:BC-inequality}
\Lambda^\ast(x) \gr (1-\varepsilon)\log \frac{1}{q_\mu(x)}+\log \frac{\varepsilon}{2^{1-\varepsilon}}
\end{equation} 
for every $x \in \operatorname{int}(\operatorname{supp}(\mu))$. Also, in~\cite[Theorem~1.1]{Pafis-2026} the second
named author obtained an analogous result for the uniform distribution on a convex body with a sharp dependence on the dimension. 

For the proof of Theorem~\ref{th:near-linear} we pass to dual radial coordinates $\theta=ru$, using $F$ as the radial variable.
We use the slicing theorem of Klartag and Lehec~\cite{Klartag-Lehec-2025}, which provides an affirmative answer to Bourgain's
slicing problem, to control the normalizing factor of this change of variables.
Theorem~\ref{th:pointwise} allows us to compare the true depth with $\widetilde q=e^{-F}/(1+J)$. Using the relative entropy between 
two tilted measures we then bound the change of a transverse covariance determinant,
producing the factor $\bigl(1+\log\lambda/(n-1)\bigr)^{n-1}$ in Theorem~\ref{th:near-linear}. 

\medskip

Two companion manuscripts use the same half-space depth wet parts for different purposes, and the results are complementary. The article~\cite{Pafis-covering} establishes half-space covering results and proves a family sensitive theorem asserting 
that entire half-spaces of small measure can be contained in the covering sets of comparable measure. In~\cite{Pafis-covering}  
it is also proved that
\begin{equation}\label{eq:intro-companion-relative}
w_\mu(s)\ls D_n\left(\frac{s}{t}\right)^{2/(n+1)}w_\mu(t), \quad 0<s\ls t<t_n.
\end{equation}
Note that applying Theorem~\ref{th:near-linear} with $\lambda=t/s$ we obtain the reverse comparison
\begin{equation}\label{eq:intro-companion-lower}
w_\mu(s)\gr C^{-n}\frac{s}{t}\left(1+\frac{\log(t/s)}{n-1}\right)^{-(n-1)}w_\mu(t), \quad 0<s<t<1.
\end{equation}

The article~\cite{Pafis-random} studies the convex hull $P_N$ of $N$ independent random points distributed according to some non-atomic probability measure $\mu$. 
It is proved that the wet-part doubling is equivalent to $E[1-\mu(P_N)]\approx w_\mu(1/N)$. Therefore, for log-concave measures, Theorem~\ref{th:near-linear} and Efron's identity give
\begin{equation}\label{eq:intro-companion-random}
\mathbb E[1-\mu(P_N)]\approx_n w_\mu(1/N),\quad\mathbb E f_0(P_N)\approx_{n}Nw_\mu(1/N),
\end{equation}
where $f_k(P_N)$ is the number of $k$-dimensional faces of $P_N$ for every $0 \ls k \ls n-1$. The results of \cite{Pafis-covering} provide the additional geometric tools for controlling the expected numbers of faces of each dimension. Combining them with Theorem~\ref{th:near-linear} yields $$\mathbb E f_k(P_N)\approx_{n,k}Nw_\mu(1/N),\quad 0\ls k\ls n-1.$$ 

\medskip

The paper is organized as follows. In Section~\ref{sec:preliminaries} we introduce the basic notation and develop the needed 
estimates for log-concave measures. In Section~\ref{sec:pointwise-estimate} we prove the depth comparison, Theorem~\ref{th:pointwise}. Section~\ref{sec:dual-radial} develops the dual radial representation. In Section~\ref{sec:near-linear} we prove the near-linear Theorem~\ref{th:near-linear}. 
Finally, in Section~\ref{sec:endpoint-counterexample} we prove that we cannot have a linear bound in Theorem~\ref{th:general-failure}.

\section{Preliminaries}\label{sec:preliminaries}

We work in ${\mathbb R}^n$, which is equipped with the standard inner product  $\langle\cdot ,\cdot\rangle $. We denote the corresponding Euclidean norm by $|\cdot |$, and write $B_2^n$ for the Euclidean unit ball, $\omega_n$ for its volume and $S^{n-1}$ for the unit sphere. Let $\sigma$ denote the surface area measure on $S^{n-1}$. For $n=1$, $S^0=\{-1,1\}$ and $\sigma$ is the counting measure. The letters $c, c_j,C,C_j$ etc. denote absolute positive constants, whose value may change from line to line. For positive real numbers $a,b$, 
the notation $a\lesssim b$ means $a\ls Cb$ for an absolute constant $C$, and $a\approx b$ means that we have both $a\lesssim b$ and $b\lesssim a$. We also write $a\approx_\xi b$ if the implicit constants depend additionally on $\xi$.

For a convex set $A \subseteq \mathbb{R}^n$ with $0$ in its interior, we define the radial function of $A$ as $\ell_A(u)=\sup \{r>0: ru \in A\} \in (0,\infty]$ for every $u \in S^{n-1}$.

A Borel probability measure $\mu$ on $\mathbb R^n$ is called log-concave if it is full-dimensional, i.e.  $\mu(H)<1$ for every hyperplane $H$ in ${\mathbb R}^n$, and
$$\mu(\lambda A+(1-\lambda)B) \gr \mu(A)^{\lambda}\mu(B)^{1-\lambda}$$
for any pair of compact sets $A,B$ in ${\mathbb R}^n$ and any $\lambda \in (0,1)$. Borell \cite{Borell-1974} has proved that, under these assumptions, $\mu $ has a log-concave density $f=e^{-V}$ with $V$ convex, which admits an exponential tail bound
$f(x)\ls ae^{-b|x|}$ for suitable $a,b>0$ (see \cite[Lemma~2.2.1]{BGVV-book}). In particular, $\mu$ has finite moments of every order. Its one-dimensional marginals are log-concave by Pr\'{e}kopa's theorem
\cite{Prekopa-1973}. Throughout the article, we work with $\mu$ log-concave. We will denote by $K$ its convex support and by $f$ its log-concave and upper-semicontinuous density representative. Thus, $f$ attains its maximum.

We say that a log-concave probability measure $\mu$ on ${\mathbb R}^n$ is centered if its barycenter ${\rm bar}(\mu)=\int_{\mathbb{R}^n}x\,d\mu(x)$ is at the origin. Moreover, $\mu$ is called symmetric if $\mu (-B)=\mu (B)$ for every Borel subset $B$ of ${\mathbb R}^n$. If $\mu$ has density $f$, we define the covariance matrix $\Cov(\mu)$ with entries
\begin{equation*}\Cov(\mu )_{ij}:=\int_{{\mathbb R}^n}x_ix_j f
(x)\,dx-\int_{{\mathbb
R}^n}x_i f (x)\,dx\int_{{\mathbb R}^n}x_j f
(x)\,dx.\end{equation*}
The isotropic constant of $\mu$ is defined by
\begin{equation*}
L_{\mu }:= \|f\|_\infty^{\frac{1}{n}} \big(\det \Cov(\mu)\big)^{\frac{1}{2n}},\end{equation*}
and $\mu $ is called isotropic if it is centered and $\Cov(\mu )=I_n$. Note that if $\mu$ is isotropic then $L_{\mu }=\|f_{\mu }\|_{\infty }^{1/n}$. 

\begin{theorem}[Klartag--Lehec]\label{th:slicing-normalization}
There exist absolute constants $0<c\ls C<\infty$ such that for every log-concave probability measure $\mu$ on $\mathbb R^n$ $$c\ls L_\mu\ls C.$$
\end{theorem}

The lower bound is simple (see \cite[Proposition~2.3.12]{BGVV-book}), while the upper bound is the affirmative answer to the slicing problem by Klartag and Lehec \cite{Klartag-Lehec-2025}.

If $\mu$ is a Borel probability measure on $\mathbb{R}^n$ and $T: \mathbb{R}^n \to \mathbb{R}^d$ is an affine function, we denote by $T_{\#}\mu$ the corresponding push-forward measure. If $\mu$ is log-concave and $T$ is surjective, then $T_{\#}\mu$ is also log-concave. 

For probability measures $\mu,\nu$, we use the relative entropy (or Kullback--Leibler divergence) 
$$D_{\mathrm{KL}}(\mu\|\nu)=\int_{\mathbb{R}^n}\log\left(\frac{d\mu}{d\nu}\right)d\mu$$
when $\mu\ll \nu$, and $D_{\mathrm{KL}}(\mu\|\nu)=+\infty$ otherwise. We note that $D_{\mathrm{KL}}(\mu\|\nu) \gr 0$ and $D_{\mathrm{KL}}(\mu\|\nu)=0$ if and only if $\mu=\nu$. 
Moreover, if $T$ is a measurable function then 
\begin{equation}\label{eq:data-processing-inequality}
D_{\mathrm{KL}}(T_{\#}\mu\|T_{\#}\nu)\ls D_{\mathrm{KL}}(\mu\|\nu).
\end{equation}
If $\mu$ has a density, we define the differential entropy of $\mu$ as $$h(\mu)=-\int f \log f.$$

\begin{theorem}\label{th:log-concave-entropy-bounds}
If $\mu$ is a log-concave probability measure on $\mathbb{R}^n$ then $$\frac12\log\bigl(\det\Cov(\mu)\bigr)-n\log C\ls h(\mu)\ls\frac12\log\bigl(\det\Cov(\mu)\bigr)+\frac n2\log(2\pi e),$$ where $C>0$ is the absolute constant of Theorem~$\ref{th:slicing-normalization}$.
\end{theorem}

For the lower bound we use the inequality $f \ls \|f\|_\infty$, as well as the upper bound of Theorem~\ref{th:slicing-normalization}, while the upper bound is the Gaussian maximum entropy inequality. For further background on information theory we refer the reader to \cite{Cover-Thomas-2006,Gray-2011}.

For $u\in S^{n-1}$ and $a\in\mathbb{R}$ let
$$H^+(u,a)=\{z:\langle u,z\rangle\gr a\}.$$
If $H^+(u,a)$ contains $x$, then $a\ls\langle u,x\rangle$, and replacing $a$ by $\langle u,x\rangle$ only decreases the measure of 
the half-space. Therefore, Tukey's half-space depth can be written as
\begin{equation}\label{eq:depth-directional}
q_\mu(x)=\inf_{u\in S^{n-1}}\mu\bigl(H^+(u,\langle u,x\rangle)\bigr),
\end{equation}
and it is enough to consider half-spaces whose boundary contains $x$.

If $\mu$ is log-concave then, for fixed $u$, the marginal $\langle u,X\rangle$ has a density, so its tail distribution function is continuous. Hence each directional half-space probability is continuous in $x$, and their infimum $q_\mu$ is upper-semicontinuous. In particular, the wet part at level $t$ $$W_\mu(t)=\{x:q_\mu(x) \ls t\}$$ is a Borel set for every $t>0$ and $w_\mu(t)=\mu(W_\mu(t))$ can be defined. If $T$ is an invertible affine transformation and $\nu=T_{\#}\mu$, then $q_\nu(x)=q_\mu(T^{-1}x)$ for every $x \in \mathbb{R}^n$. Therefore, $w_\mu(t)=w_\nu(t)$ for every $t>0$.

We will also use the bound
\begin{equation}\label{eq:w-lower}
w_\mu(t)\gr t, \quad 0<t<1,
\end{equation}
for every log-concave probability measure. For a proof of~\eqref{eq:w-lower}, choose a half-space whose measure
is exactly equal to $t$. Then, for every point $x$ of that half-space we have $q_{\mu}(x)\ls t$. Such a half-space exists 
because the tail distribution function of every marginal is continuous and ranges from $0$ to $1$.

\smallskip 

Given a convex function $g:{\mathbb R}^n\to (-\infty ,\infty ]$, the Legendre--Fenchel transform $g^\ast$
of $g$ is defined by
\begin{equation*}g^\ast(x):=\sup_{\theta}\{\langle\theta,x\rangle-g(\theta)\}\end{equation*} and is also a convex function.
Let now $\mu$ be a probability measure on $\mathbb{R}^n$. Its log-Laplace transform is defined by
$$\Lambda(\theta)=\log\int_{\mathbb{R}^n} e^{\langle\theta,x\rangle}\,d\mu(x)$$
on its effective domain $\Omega=\{\theta: \Lambda(\theta)<\infty\}$. The Cram\'er transform is the Legendre--Fenchel transform of $\Lambda$, defined by
$$\Lambda^*(x)=\sup_{\theta}\{\langle\theta,x\rangle-\Lambda(\theta)\}.$$ 
Both $\Lambda$ and $\Lambda^\ast$ are convex and lower-semicontinuous functions. Standard references for Legendre duality and exponential families are \cite{Rockafellar-1970,Brown-1986}. For every $\theta \in \Omega$, we define the exponential tilt
\begin{equation}\label{eq:tilt-def}
d\mu_\theta(z)=e^{\langle\theta,z\rangle-\Lambda(\theta)}\,d\mu(z).
\end{equation}
In the next proposition we gather some useful properties of $\Lambda$ and $\Lambda^*$.

\begin{proposition}\label{prop:lambda-cramer-properties}
Let $\mu$ be a log-concave probability measure on $\mathbb{R}^n$ with $K={\rm supp}(\mu)$. Then:
\begin{enumerate}
\item[{\rm (i)}] The set $\Omega$ is open and convex, and $\Lambda$ is $C^{\infty}$ and strictly convex on $\Omega$, with positive definite Hessian.
\item[{\rm (ii)}] The function $\Lambda^*: {\rm int}(K) \to \mathbb{R}$ is $C^{\infty}$.
\item[{\rm (iii)}] The map $\nabla \Lambda: \Omega \to {\rm int}(K)$ is a smooth diffeomorphism from $\Omega$ onto ${\rm int}(K)$ and $(\nabla \Lambda)^{-1}=\nabla \Lambda^*$.
\item[{\rm (iv)}] For every $x \in {\rm int}(K)$ we have $\Lambda^*(x)=\langle x,\theta\rangle -\Lambda(\theta)$ if and only if $x=\nabla\Lambda(\theta)$.
\item[{\rm (v)}] For every $\theta$ in $\Omega$:
$$\mathbb{E}_{\mu_\theta}(X)=\nabla\Lambda(\theta),\quad \Cov_{\mu_\theta}(X)=\nabla^2\Lambda(\theta)=:\Hess(\theta).$$
\end{enumerate}
\end{proposition}

\begin{proof}The effective domain $\Omega$ is convex by H\"{o}lder's inequality and contains $0$. It is also open: for any $\theta_0\in\Omega$, 
the tilted distribution $\mu_{\theta_0}$ is again a log-concave probability measure, hence has a density bounded by $a e^{-b|z|}$ for some $a,b>0$. 
It follows that $\mathbb E_{\mu_{\theta_0}}e^{\langle h,X\rangle}<\infty$ whenever $|h|<b$, and then $\Lambda(\theta_0+h)<\infty$. Therefore, 
on $\Omega$ we can differentiate under the integral, and this proves smoothness and the identities in ${\rm (v)}$.

The covariance matrix is positive definite, because $\mu_\theta$ is full-dimensional, which means that it is not supported on a hyperplane. Thus $\Lambda$ is strictly convex.

We prove surjectivity directly. Fix $x\in{\rm int}(K)$ and consider the  function
$\phi_x(\theta)=\Lambda(\theta)-\langle\theta,x\rangle$ on all of $\mathbb R^n$, taking the value $+\infty$ outside $\Omega$.
It is proper and lower-semicontinuous. Choose $r>0$ so that $\overline B(x,2r)\subseteq {\rm int}(K)$. If $f$ is the log-concave density of $\mu$, then $f$ is strictly positive and continuous on
${\rm int}(K)$. In particular, there exists $m>0$ such that $f\gr m$ on $\overline B(x,2r)$. For every $\theta\ne0$, we set $e=\theta/|\theta|$ and integrating over $B(x+re,r/2)$ we obtain
\begin{equation}\label{eq:log-laplace-coercivity}
\phi_x(\theta)=\log\int_{\mathbb{R}^n} e^{\langle\theta,z-x\rangle}f(z)\,dz\gr \frac r2|\theta|+\log\bigl(m\,{\rm vol}_n(\tfrac{r}{2}B_2^n)\bigr).
\end{equation}
This shows that $\phi_x$ is coercive on $\mathbb R^n$ and, as it is lower-semicontinuous, $\phi_x$ attains its minimum. Since $\phi_x(0)=0$, every minimizer has finite value and therefore belongs to the open set $\Omega$. Strict convexity implies that this minimizer $\theta_x$ is unique, 
and differentiating at this point we see that $\nabla\Lambda(\theta_x)=x$.

Conversely, $\nabla\Lambda(\theta)$ belongs to
${\rm int}(K)$ for every $\theta \in \Omega$, since it is equal to the barycenter of $\mu_\theta$, which is a log-concave probability measure supported on $K$. The inverse function theorem and the positive definite Hessian now show that $\nabla\Lambda$ is a smooth diffeomorphism onto ${\rm int}(K)$. The inverse is $\nabla\Lambda^*$ by the Fenchel identity, which also proves ${\rm (ii)}$.
\end{proof}

Related general convex analytic formulations can be found in \cite{Brown-1986,Rockafellar-1970}. For $\theta \in \Omega$ and $v \in \mathbb{R}^n$ we write
\begin{equation}\label{eq:localnorm}
\|v\|_\theta^2=v^T\Hess(\theta)v,\quad J(\theta)=\|\theta\|_\theta.
\end{equation} Observe that $$\|v\|_\theta^2=\Var_{\mu_\theta}\langle v,X\rangle.$$
Next we list some one-dimensional estimates that will be used below.

\begin{lemma}\label{lem:1d-density}
Let $Z$ be a centered one-dimensional log-concave random variable with variance $\sigma^2>0$ and density $g$. Then there exist 
absolute constants $c,C,\delta>0$ such that
\begin{equation}\label{eq:1d-density-bounds}
\|g\|_\infty\ls \frac C\sigma\quad\text{and}\quad g(s)\gr\frac c\sigma \quad\text{for }0\ls s\ls\delta\sigma.
\end{equation}
Therefore,
\begin{equation}\label{eq:saddlepoint-1d}
\mathbb{E}[e^{-Z}\mathds{1}_{\{Z\gr0\}}]\approx\frac1{1+\sigma}.
\end{equation}
\end{lemma}

\begin{proof}After scaling, we may assume that $\sigma=1$. For an isotropic density $g$ on
$\mathbb R$, Theorem~5.14(a),(d),(e) of Lov\'{a}sz and Vempala \cite{Lovasz-Vempala-2007} gives the following inequalities:
$$g(t)\gr 2^{-9|t|}g(0)\quad (|t|\ls1/9),\quad g(0)\gr2^{-7},\quad \|g\|_\infty\ls2^8.$$
In particular, $g(t)\gr2^{-8}$ for $0\ls t\ls1/9$. For general variance $\sigma^2$, rescaling we get 
$$g(t)\gr\frac{2^{-8}}\sigma\quad (0\ls t\ls\sigma/9),\quad \|g\|_\infty\ls\frac{2^8}{\sigma}.$$
Thus \eqref{eq:1d-density-bounds} holds with $c=2^{-8}$, $C=2^8$ and $\delta=1/9$.

Next we prove \eqref{eq:saddlepoint-1d}. For the upper bound we write
$$\mathbb{E}[e^{-Z}\mathds{1}_{\{Z\gr0\}}]\ls \min\left\{1,\int_0^\infty e^{-s}\|g\|_\infty\,ds\right\}\lesssim \min\{1,\sigma^{-1}\}.$$
For the lower bound, we integrate on $[0,\eta\min\{1,\sigma\}]$ choosing $\eta>0$ so that this interval lies in $[0,\delta\sigma]$. 
On this interval we have $e^{-s}\gtrsim1$ and $g(s)\gtrsim\sigma^{-1}$, which gives
$$\mathbb{E}[e^{-Z}\mathds{1}_{\{Z\gr0\}}]\gtrsim \min\{1,\sigma^{-1}\}.$$
Since $\min\{1,\sigma^{-1}\}\approx(1+\sigma)^{-1}$, the result follows.
\end{proof}

We shall use the following form of Gr\"unbaum's inequality: if $Z$ is a one-dimensional log-concave random variable, then
\begin{equation}\label{eq:grunbaum}
\mathbb{P}(Z\gr\mathbb{E} Z)\gr e^{-1},\quad\mathbb{P}(Z\ls\mathbb{E} Z)\gr e^{-1}.
\end{equation}
This is the one-dimensional case of Gr\"unbaum's centroid half-space inequality \cite{Grunbaum-1960}; see also \cite{Lovasz-Vempala-2007}.

We shall also use the following third moment estimate (third derivative estimates of this type appear in self-concordance theory; see \cite{Bubeck-Eldan-2019,Lee-Yue-2021}).

\begin{lemma}\label{lem:skew}
Let $Z$ be a log-concave, one-dimensional random variable. Then,
\begin{equation}\label{eq:skew}
\big|\mathbb{E}(Z-\mathbb{E} Z)^3\big|\ls2\,\Var(Z)^{3/2}.
\end{equation}
\end{lemma}

\begin{proof}Translate $Z$ so that $\mathbb{E} Z=0$. From Lemma~2 of Bubeck and Eldan \cite{Bubeck-Eldan-2019} we have
$$\mathbb{E} Z^3\ls2(\mathbb{E} Z^2)^{3/2}.$$
Applying the same inequality to $-Z$ we obtain the lower bound. This proves \eqref{eq:skew}.
\end{proof}

Finally, we will need a relative entropy comparison of one-dimensional variances.

\begin{lemma}\label{lem:var-kl}Let $\mu,\nu$ be log-concave probability measures on $\mathbb{R}^n$ with $D_{\mathrm{KL}}(\mu\|\nu)<\infty$, and let
$v\in\mathbb{R}^n$. Then,
\begin{equation}\label{eq:var-kl}
\sqrt{\Var_\nu\langle v,X\rangle} \ls C e^{D_{\mathrm{KL}}(\mu\|\nu)} \sqrt{\Var_\mu\langle v,X\rangle},
\end{equation}
where $C>0$ is an absolute constant.
\end{lemma}

\begin{proof}The case $v=0$ is trivial, so we may assume that $v \neq 0$. Set $D=D_{\mathrm{KL}}(\mu\|\nu)<\infty$, and let ${T_{v}}_{\#}\mu,{T_{v}}_{\#}\nu$ be the
push-forwards under $T_v(x)=\langle v,x\rangle$. Since $\mu$ and $\nu$ are full-dimensional, both projections have strictly positive
variance. From \eqref{eq:data-processing-inequality} we obtain
$$D_{\mathrm{KL}}({T_{v}}_{\#}\mu\|{T_{v}}_{\#}\nu)\ls D.$$
The projections ${T_{v}}_{\#}\mu,{T_{v}}_{\#}\nu$ have one-dimensional log-concave densities $p_v,q_v$. We write their standard deviations as $\sigma_{p_v},\sigma_{q_v}$. Using Theorem~\ref{th:log-concave-entropy-bounds} in dimension $1$ we get
$$h({T_{v}}_{\#}\mu)\ls \frac12\log(2\pi e\sigma_{p_v}^2).$$ Applying Lemma~\ref{lem:1d-density} to the centered translation of $q_v$ we deduce that
$$\|q_v\|_\infty\ls C_1\sigma_{q_v}^{-1}.$$ Since $\log q\ls\log\|q\|_\infty$, we get
\begin{align*}
D_{\mathrm{KL}}({T_{v}}_{\#}\mu\|{T_{v}}_{\#}\nu) &=-h({T_{v}}_{\#}\mu)-\int p_v\log q_v\gr-\frac12\log(2\pi e\sigma_{p_v}^2)-\log\|q_v\|_\infty\\
&\gr\log\frac{\sigma_{q_v}}{\sigma_{p_v}}-C_1.
\end{align*}
Therefore, $\sigma_{q_v}\ls e^{C_1}e^D\sigma_{p_v}$. This proves \eqref{eq:var-kl}.
\end{proof}

\section{Half-space depth}\label{sec:pointwise-estimate}

We keep the notation of Section~\ref{sec:preliminaries}. In particular, $x(\theta)$ and $\Hess(\theta)$ are the mean and covariance of $\mu_\theta$.

\subsection{Half-space probabilities}\label{sec:pointwise}

Fix $u\in S^{n-1}$, $x\in {\rm int}(K)$, and set $a=\langle u,x\rangle$. Let
$$H=H^+(u,a)=\{z:\langle u,z\rangle\gr a\}.$$
Define the log-Laplace transform of the one-dimensional marginal $Y=\langle u,X\rangle$ by
$$\psi(s)=\Lambda(su).$$
Since $x$ lies in the interior of the full-dimensional support $K$, there exists $\delta>0$ such that $x+su\in K$ for all $|s|<\delta$. It follows that the support interval of $Y$ contains $(a-\delta,a+\delta)$, so $a$ belongs to its interior. Applying the one-dimensional form of Proposition~\ref{prop:lambda-cramer-properties}, we see that there exists a unique $\tau\in{\rm dom}(\psi)$ such that
\begin{equation}\label{eq:1d-saddle-eq}
\psi'(\tau)=a.
\end{equation}
Suppose now that $\mu(H)<e^{-1}$. By \eqref{eq:grunbaum}, we deduce that $a>\mathbb{E}Y$. Since $\psi$ is strictly convex, we
have that $\psi'$ is strictly increasing. Thus, since $\psi'(0)=\mathbb{E}Y$, the parameter $\tau$ defined by \eqref{eq:1d-saddle-eq} is positive. 
We set
$$\zeta=\tau u,\quad I_\zeta=\tau a-\Lambda(\zeta), \quad J(\zeta)^2=\zeta^T\Hess(\zeta)\zeta.$$
Note that $\mathbb E_{\mu_\zeta}Y=\psi'(\tau)=a$, and hence $S=\tau(Y-a)$ is centered.  It is
also log-concave and its variance equals
$$\Var_{\mu_\zeta}(S)=\tau^2\psi''(\tau)=J(\zeta)^2.$$
Since 
$d\mu(x)=e^{-\langle\zeta,x\rangle+\Lambda(\zeta)}d\mu_\zeta(x)$ we get
\begin{align}\label{eq:halfspace-tilt-identity}
\mu(H)=e^{-I_\zeta}\mathbb E_{\mu_\zeta}[e^{-S}\mathds{1}_{\{S\gr0\}}].
\end{align}
Applying Lemma~\ref{lem:1d-density} we obtain
\begin{equation}\label{eq:halfspace-saddle}
\mu(H)\approx \frac{e^{-I_\zeta}}{1+J(\zeta)}.
\end{equation}

The next lemma is used to compare exponential tilts. Its role is to convert a relative entropy cost into a bound in the
covariance norm at the endpoint of a segment. The essential observation is that the reciprocal square root of the second derivative is
$1$-Lipschitz.

\begin{lemma}\label{lem:increment}
Let $g:[0,1]\to\mathbb{R}$ be $C^3$, with $g(0)=g'(0)=0$ and $g''>0$. Assume that
\begin{equation}\label{eq:sc-line}
|g'''(s)|\ls2g''(s)^{3/2} \quad(0\ls s\ls1).
\end{equation}
If $D=g(1)$ and $\gamma=\sqrt{g''(1)}$, then
\begin{equation}\label{eq:increment-sharp}
D\gr\log(1+\gamma)-\frac{\gamma}{1+\gamma},
\end{equation}
and in particular
\begin{equation}\label{eq:A-expD}
\gamma\ls e^{D+1}.
\end{equation}
\end{lemma}

\begin{proof}Set $a(s)=\sqrt{g''(s)}$. Differentiating $a^{-1}$ and using \eqref{eq:sc-line} we get
$$\left|\frac d{ds}\frac1{a(s)}\right|=\frac{|g'''(s)|}{2g''(s)^{3/2}}\ls1.$$
Since $\gamma=a(1)$, integrating the inequality from $s$ to $1$ gives
$$\frac1{a(s)}\ls \frac1\gamma+1-s.$$ Therefore, as $g''(s)=a(s)^2$ we have $$g''(s) \gr \frac{\gamma^2}{(1+\gamma(1-s))^2}.$$
Using $g(0)=g'(0)=0$ and integrating by parts, we obtain
$$D=g(1)=\int_0^1(1-s)g''(s)\,ds\gr \int_0^1\frac{(1-s)\gamma^2}{(1+\gamma(1-s))^2}\,ds=\log(1+\gamma)-\frac{\gamma}{1+\gamma}.$$
Since $\gamma/(1+\gamma)\ls1$, \eqref{eq:A-expD} follows.
\end{proof}

\subsection{Comparison of directions}

In this subsection we prove Theorem~\ref{th:pointwise}. The upper bound uses one explicit half-space. For the lower bound, we compare an arbitrary half-space tilt with the tilt whose mean is
$x$. The rate advantage of the former is exactly a relative entropy $D$. The next argument shows that its variance factor costs at most
$e^D$. These two factors cancel.

\begin{proof}[Proof of Theorem~$\ref{th:pointwise}$]
We consider first the case $\theta=0$. Then $x=\mathbb{E} X$. From Gr\"unbaum's inequality we have $q_\mu(x)\gr e^{-1}$. On the
other hand, at least one of the two opposite half-spaces through $x$ has measure at most $1/2$. Since $F(0)=J(0)=0$, \eqref{eq:pointwise} follows. 

From now on we assume that $\theta\ne 0$. By \eqref{eq:depth-directional}, we may consider only half-spaces through $x$.

\medskip

\noindent \textbf{Upper bound.} Consider the half-space
$$H_\theta=\{z:\langle\theta,z\rangle\gr\langle\theta,x\rangle\}.$$
The random variable $S_\theta=\langle\theta,X-x\rangle$ is centered and has variance $J(\theta)^2$
with respect to $\mu_\theta$. It follows that
$$\mu(H_\theta)=e^{-F(\theta)} \mathbb{E}_{\mu_\theta}[e^{-S_\theta}\mathds{1}_{\{S_\theta\gr0\}}].$$
Therefore, using Lemma~\ref{lem:1d-density} we obtain
$$q_\mu(x)\ls\mu(H_\theta)\lesssim\frac{e^{-F(\theta)}}{1+J(\theta)}.$$

\medskip

\noindent \textbf{Lower bound.} Let $H=H^+(u,\langle u,x\rangle)$ be any half-space through $x$. If $\mu(H)\gr e^{-1}$, then the lower bound 
follows after choosing the implicit absolute constant small enough, because $F\gr0$ and $(1+J(\theta))^{-1}\ls1$. We may therefore assume that
$\mu(H)<e^{-1}$. We construct its saddlepoint tilt $\zeta=\tau u$ as in \eqref{eq:1d-saddle-eq}. By \eqref{eq:halfspace-saddle},
\begin{equation}\label{eq:competing-saddle}
\mu(H)\approx\frac{e^{-I_\zeta}}{1+J(\zeta)}.
\end{equation}
Both $\mu_\theta$ and $\mu_\zeta$ have positive densities with respect to $\mu$, hence they are mutually absolutely continuous. Since
$$\log\frac{d\mu_\theta}{d\mu_\zeta}(z)=\langle z,\theta-\zeta\rangle-\Lambda(\theta)+\Lambda(\zeta)$$
and $\mu_\theta$ has mean $x$, direct integration shows that
\begin{equation}\label{eq:KL-rate-gap}
D:=D_{\mathrm{KL}}(\mu_\theta\|\mu_\zeta) =\langle x,\theta-\zeta\rangle-\Lambda(\theta)+\Lambda(\zeta) =F(\theta)-I_\zeta.
\end{equation}
We next control $J(\zeta)=\|\zeta\|_\zeta$. Applying Lemma~\ref{lem:var-kl} in direction $\theta$ we get
\begin{equation}\label{eq:theta-under-Q}
\|\theta\|_\zeta \ls C e^D\|\theta\|_\theta =C e^D J(\theta).
\end{equation}
Let $h=\zeta-\theta$ and observe that $\theta+sh \in \Omega$ for every $s \in [0,1]$. If $h=0$, then \eqref{eq:h-under-Q} below is trivial. 
So, we assume that $h\neq0$ and define
$$g(s)=\Lambda(\theta+sh)-\Lambda(\theta)-s\langle h,x\rangle,\quad 0\ls s\ls1.$$
Then $g(0)=g'(0)=0$, $g(1)=D$, and
$$g''(s)=h^T\Hess(\theta+sh)h.$$
Moreover $g'''(s)$ is the third centered moment of $\langle h,X\rangle$ under the tilt $\mu_{\theta+sh}$. The tilted distribution remains 
log-concave, hence Lemma~\ref{lem:skew} gives
$$|g'''(s)|\ls2g''(s)^{3/2}.$$
Applying Lemma~\ref{lem:increment} we obtain
\begin{equation}\label{eq:h-under-Q}
\|h\|_\zeta=\sqrt{g''(1)}\ls e^{D+1}.
\end{equation}
Combining \eqref{eq:theta-under-Q} and \eqref{eq:h-under-Q} we see that
\begin{equation}\label{eq:Jcomp}
1+J(\zeta)=1+\|\zeta\|_\zeta\ls1+\|\theta\|_\zeta+\|h\|_\zeta\ls C e^D(1+J(\theta)),
\end{equation}
because $D\gr0$. Finally, from \eqref{eq:competing-saddle}, \eqref{eq:KL-rate-gap} and \eqref{eq:Jcomp} we see that
$$\mu(H)\gtrsim\frac{e^{-F(\theta)+D}}{e^D(1+J(\theta))}=\frac{e^{-F(\theta)}}{1+J(\theta)}.$$
The bound is uniform over all half-spaces $H$ through $x$. Taking the infimum over $u$ we complete the proof of the theorem.
\end{proof}

\begin{remark}\rm
The next example illustrates the need for a correction factor. Let $\mu$ be the one-dimensional Laplace probability measure with density $\frac12e^{-|y|}$. For every $r\gr0$,
$$q_\mu(r)=\frac12e^{-r},\quad \Lambda(\theta)=-\log(1-\theta^2).$$
The non-negative saddlepoint $\theta_r$ satisfies the equation
$$r=\Lambda'(\theta_r)=\frac{2\theta_r}{1-\theta_r^2},$$
so
$$\theta_r=\frac{r}{1+\sqrt{1+r^2}},\quad r\gr0.$$
Note that
$$F(\theta_r)=r-\log r+O(1)\quad\text{and}\quad J(\theta_r)=r+O(1)$$
as $r\to\infty$. Hence
$$\frac{e^{-F(\theta_r)}}{1+J(\theta_r)}\approx e^{-r}\approx q_\mu(r).$$
Here $q_\mu(r)e^{F(\theta_r)}\approx r^{-1}$. Note that a factor depending only on the dimension in place of $1+J$ cannot give a two-sided 
pointwise estimate for this example.
\end{remark}

\section{Dual coordinates and radial regularity}\label{sec:dual-radial}

We now pass from the pointwise depth estimate to a comparison of wet-part measures. The main task is to control how the measure changes along a
ray in the tilt domain $\Omega$.

\subsection{Dual coordinates}\label{sec:dual-polar}

Our aim is to express $\mu$ in tilt coordinates, then in polar coordinates, and finally in the rate coordinate $F$. We first check that
$F$ takes all values in $(0, \infty)$ along every ray. The next lemma describes the behavior of $\Lambda^*$ outside compact subsets of ${\rm int}(K)$.

\begin{lemma}\label{lem:boundary-star-explosion}
We have $\Lambda^*(x)\to\infty$ as $|x|\to\infty$ and $\Lambda^\ast(x)=\infty$ for every $x \in \mathbb{R}^n \setminus {\rm int}(K)$. In particular, for every $u\in S^{n-1}$,
$$\lim_{r \to \ell_{\Omega}(u)^-} \Lambda^*(\nabla \Lambda(ru))=\infty.$$
\end{lemma}

\begin{proof}
For the first assertion, since $\Omega$ is open and contains $0$, by the continuity of $\Lambda$, there exist $m,M>0$ such that $\Lambda(\theta)\ls M$ for every $\theta \in mB_2^n$. Let $x\neq 0$ and set $\theta =mx/|x|$. Then, 
$$\Lambda^*(x) \gr \langle x,\theta\rangle -\Lambda(\theta)\gr m|x|-M.$$
Taking $|x| \to \infty$ we obtain the result. 
    
For the second assertion, let $x \in \mathbb{R}^n \setminus {\rm int}(K)$. Then we can find some $v \in S^{n-1}$ such that $\langle v,z\rangle \ls \langle v,x\rangle$  for every $z \in K$. Since $\mu$ is absolutely continuous and $v\neq 0$, we have $\mu(\{z:\langle v,z\rangle =\langle v,x\rangle\})=0$. Hence, by the dominated convergence theorem we obtain $$\exp(\Lambda (tv)-t\langle v,x\rangle)=\mathbb{E}\,\left[e^{t\langle v,X-x\rangle}\right] \to 0,$$ 
as $t \to \infty$. Therefore, 
$$\langle x,tv\rangle -\Lambda(tv) \to +\infty,$$ 
as $t \to \infty$, which implies that $\Lambda^\ast(x)=+\infty$.

To prove the final assertion, observe that for every $M>0$ the level set 
$$B_M(\mu)=\{x:\Lambda^\ast(x) \ls M\}$$ is a compact subset of ${\rm int}(K)$ as it is closed by lower semicontinuity and bounded by the first part of
the proof. From the second part we see that it is contained in ${\rm int}(K)$.

The last claim follows from the fact that $ru$ eventually escapes every compact subset of ${\rm dom}(\Lambda)$, as $r \to \ell_{{\rm dom}(\Lambda)}(u)^-$, combined with Proposition~\ref{prop:lambda-cramer-properties} {\rm (iii)}.
\end{proof}

We remind now that for some $\theta \in \Omega$ the tilted density is
$$f_\theta(x)=e^{\langle\theta,x\rangle-\Lambda(\theta)}f(x).$$
At the barycenter $x(\theta)$ of $\mu_\theta$ we have
\begin{equation}\label{eq:tiltdensity-bary}
f_\theta(x(\theta))=e^{F(\theta)}f(x(\theta)).
\end{equation}
Define
\begin{equation}\label{eq:b-factor}
b(\theta)=f_\theta(x(\theta))\sqrt{\det\Hess(\theta)}.
\end{equation}

\begin{proposition}\label{prop:dual-polar}
For every $\theta\in \Omega$,
\begin{equation}\label{eq:dual-density-exact}
d\mu(x(\theta))=e^{-F(\theta)}b(\theta)\sqrt{\det\Hess(\theta)}\,d\theta.
\end{equation}
For $\theta=ru$, $r>0$, $u\in S^{n-1}$, set
\begin{equation}\label{eq:aJ}
a(r,u)=u^T\Hess(ru)u\quad\text{and}\quad J(r,u)=r\sqrt{a(r,u)}.
\end{equation}
Then
\begin{equation}\label{eq:dFdr}
\frac{d}{dr}F(ru)=r a(r,u).
\end{equation}
For every $n\gr1$ define
\begin{equation}\label{eq:Kn}
\kappa_n(r,u)=\frac{\sqrt{\det\Hess(ru)}\,r^{n-2}}{a(r,u)},
\end{equation}
which reduces to $\kappa_1(r,u)=1/J(r,u)$ when $n=1$. For every non-negative Borel measurable function $\Phi:\mathbb{R}^n\to\mathbb{R}$ we have
\begin{equation}\label{eq:dual-polar-disintegration}
\int_{\mathbb{R}^n}\Phi(x)\,d\mu(x)=\int_{S^{n-1}}\int_0^\infty\Phi(x(r(F,u)u))p_u(F)\,dF\,d\sigma(u),
\end{equation}
where $r=r(F,u)$ is the inverse of $r\mapsto F(ru)$ and  
\begin{equation}\label{eq:ray-density-entropic}
p_u(F)=b(r(F,u)u)\rho_u(F),\quad
\rho_u(F)=e^{-F}\kappa_n(r(F,u),u).
\end{equation}
In particular, we have
\begin{equation}\label{eq:ray-density-comparison}
c_1^n\rho_u(F)\ls p_u(F)\ls c_2^n\rho_u(F),
\end{equation}
where $0<c_1\ls c_2<\infty$ are  absolute constants.
\end{proposition}

\begin{proof}Since $x(\theta)=\nabla\Lambda(\theta)$, its derivative is $\Hess(\theta)$ and its Jacobian is
positive. This shows that
$$\int_{\mathbb{R}^n}\Phi(x)\,d\mu(x)=\int_{\Omega}\Phi(x(\theta))f(x(\theta))\det\Hess(\theta)\,d\theta.$$
From \eqref{eq:tiltdensity-bary} and the definition of $b$ in \eqref{eq:b-factor} we see that
$$f(x(\theta))\det\Hess(\theta)=e^{-F(\theta)}\,\bigl[f_\theta(x(\theta))\sqrt{\det\Hess(\theta)}\bigr]\sqrt{\det\Hess(\theta)}.$$
This proves \eqref{eq:dual-density-exact}.

For any fixed direction $u$ we have
$$F(ru)=r\langle u,\nabla\Lambda(ru)\rangle-\Lambda(ru), \quad 0<r<\ell_{\Omega}(u),$$
and differentiating we obtain
$$\frac{dF(ru)}{dr} =r\,a(r,u)>0.$$ 
Since $a(r,u)>0$, the function $r\mapsto F(ru)$ is strictly increasing. By Lemma~\ref{lem:boundary-star-explosion}, it maps
$(0,\ell_{\Omega}(u))$ onto $(0,\infty)$. It follows that the inverse $r(F,u)$ exists for each $F>0$ and has derivative $dr/dF=[ra(r,u)]^{-1}$. 

Polar coordinates give $d\theta=r^{n-1}dr\,d\sigma(u)$, while $$dr\,d\sigma(u)=\frac{dF\,d\sigma(u)}{ra(r,u)}.$$
Substituting into \eqref{eq:dual-density-exact} we get \eqref{eq:dual-polar-disintegration}.
In dimension one the polar formula is the sum over $u\in\{-1,1\}$. Hence, $\det\Hess(ru)=a(r,u)$ and $r^{n-2}=r^{-1}$ which imply that $\kappa_1(r,u)=1/J(r,u)$.

Finally, since ${\rm bar}(\mu_\theta)=x(\theta)$, Fradelizi's inequality (see \cite[Theorem~2.2.2]{BGVV-book}) gives $$f_\theta(x_\theta) \ls \|f_\theta\|_\infty \ls e^n f_\theta(x_\theta).$$ Therefore, as ${\rm Cov}(\mu_\theta)=\Hess(\theta)$, we can write $$e^{-n}L_{\mu_\theta}^n\ls b(\theta) \ls L_{\mu_\theta}^n.$$
Applying Theorem~\ref{th:slicing-normalization} we obtain~\eqref{eq:ray-density-comparison}.\end{proof}

We define
\begin{equation}\label{eq:qtilde}
\widetilde q(x(ru))=\frac{e^{-F(ru)}}{1+J(r,u)}, \quad u \in S^{n-1},\ 0<r <\ell_\Omega(u).
\end{equation}
At the barycenter, we define $\widetilde q(x(0))=1$. This defines $\widetilde q$ on ${\rm int}(K)$ by the dual diffeomorphism.
We extend $\widetilde q$ to $\mathbb R^n\setminus{\rm int}(K)$ by setting $\widetilde q(x)=0$. The actual half-space depth is also zero there:
for a point outside $K$ this follows by a separation argument, while at a boundary point a supporting half-space has $\mu$-measure zero, since $\mu$ is
absolutely continuous. Therefore, Theorem~\ref{th:pointwise} shows that $c_0\widetilde q\ls q_\mu\ls C_0\widetilde q$,
with the same absolute constants. For every $t>0$ we write
\begin{equation}\label{eq:modified-wet-def}
\widetilde w_\mu(t)=\mu(\{x:\widetilde q(x)\ls t\}).
\end{equation} We also set
\begin{equation}\label{eq:Gdef}
G(r,u)=-\log\widetilde q(x(ru))=F(ru)+\log(1+J(r,u)).
\end{equation}

\subsection{Radial derivatives}\label{sec:ray-analysis}

We shall use two facts: $G$ increases along every ray, and if $G\gr 1$ then $dG/dF$ is bounded above and below.

We fix $u$ and omit it when there is no ambiguity. We write $a_r=da/dr$ for the derivative with respect to the original radius, while
$J'=dJ/dF$ and $G'=dG/dF$ denote derivatives with respect to the rate coordinate. From \eqref{eq:dFdr} we see that
$$\frac d{dF}=\frac1{ra}\frac d{dr}.$$ 
Differentiating $a(r)=u^T\Hess(ru)u$, we get
$$a_r(r)=D^3\Lambda(ru)[u,u,u].$$
This is the third centered moment of $\langle u,X\rangle$ with respect to $\mu_{ru}$, while $a(r)$ is its variance. 
Then, from Lemma~\ref{lem:skew} we get
\begin{equation}\label{eq:a-prime}
|a_r(r)|\ls2a(r)^{3/2}.
\end{equation}
We write
\begin{equation}\label{eq:epsilon}
\varepsilon(r,u)=\frac{a_r(r,u)}{2a(r,u)^{3/2}},\quad |\varepsilon|\ls1.
\end{equation}
Since $J=r\sqrt a$ and $dF/dr=ra$, we see that
\begin{equation}\label{eq:dlogJdF}
\frac{d\log J}{dF}=\frac1{J^2}+\frac{\varepsilon}{J},
\end{equation}
Indeed, $\frac{d}{dr}\log J=r^{-1}+a_r/(2a)$. Dividing by $ra$ we obtain \eqref{eq:dlogJdF} because $r^2a=J^2$ and
$a_r=2\varepsilon a^{3/2}$. It follows that $J'=J^{-1}+\varepsilon$, and hence
\begin{equation}\label{eq:dGdF}
G'=1+\frac{1/J+\varepsilon}{1+J}.
\end{equation}
Since $\varepsilon\gr-1$, we have
\begin{equation}\label{eq:Glower}
G'\gr\frac{J^2+1}{J(J+1)}\gr2(\sqrt2-1)=:c_0>0.
\end{equation}
In particular, $G$ is strictly increasing along every dual ray.

Also $G\to0$ as $r \to 0$, while $G\gr F\to\infty$
at the other endpoint. Hence $G$ maps each ray onto $(0,\infty)$, so its inverse is defined at every positive level.

If $J\gr1/2$, then as $\varepsilon \ls 1$ we obtain
\begin{equation}\label{eq:Gupper}
G'\ls1+\frac{1/J+1}{1+J}=1+\frac1J\ls3.
\end{equation}
We show that $G\gr1$ implies $J\gr1/2$. Since
$$\frac{dJ}{d\log r}=J(1+\varepsilon J),$$
$J$ is strictly increasing as long as $J<1$. This shows that the level $1/2$
can only be crossed upwards as $r$ increases: if $J(r)\ls1/2$ at some $r>0$, then $J(s)\ls1/2$ for all $0<s\ls r$. Moreover,
$$\frac{dF}{d\log r}=J^2.$$
Therefore, when $J\ls1/2$,
\begin{equation}\label{eq:dFdJ}
\frac{dF}{dJ}=\frac{J}{1+\varepsilon J}\ls\frac{J}{1-J}.
\end{equation}
As $r\to0$, both $F$ and $J$ tend to zero. Integrating \eqref{eq:dFdJ} we get
\begin{equation}\label{eq:FsmallJ}
F\ls -J-\log(1-J)\quad(J\ls1/2).
\end{equation}
Consequently,
$$G=F+\log(1+J)\ls -J-\log(1-J)+\log(1+J).$$
For $J \in [0,1/2]$ the right-hand side is an increasing function and at $J=1/2$ it becomes $-1/2+\log3<1$. Hence
\begin{equation}\label{eq:GimpliesJ}
G\gr1\quad\Longrightarrow\quad J\gr\frac12.
\end{equation}
Combining \eqref{eq:Glower}, \eqref{eq:Gupper} and \eqref{eq:GimpliesJ}, on the entire tail $G\gr1$ we have
\begin{equation*}
c_0\ls G'\ls3.
\end{equation*}
We summarize the main results of this subsection into the following lemma.

\begin{lemma}\label{lem:G-main-properties}
Fix $u \in S^{n-1}$. We have that $G$, as a function of $F$, is strictly increasing and maps $(0,\infty)$ to $(0,\infty)$. 
Moreover, on $G \gr 1$ $$c_0 \ls G' \ls 3, \quad J \gr \frac{1}{2},$$ where $c_0>0$ is an absolute constant.
\end{lemma}

\section{Near-linear doubling}\label{sec:near-linear}

We prove Theorem~\ref{th:near-linear} by comparing two exponential tilts on the same dual ray. To do so, we first need to appropriately control the ratio of the determinants of their covariance matrices.

\subsection{Controlling covariance determinants via relative entropy}
\begin{lemma}\label{lem:entropy-volume-comparison}
If $\mu,\nu$ are log-concave probability measures on $\mathbb R^n$ with $D=D(\mu\Vert\nu)<\infty$, then
\begin{equation}\label{eq:entropy-difference-bound}
 h(\mu)-h(\nu)\ls2n+n\log(1+D/n).
\end{equation}
Consequently, there exists an absolute constant $C>0$ such that
\begin{equation}\label{eq:covariance-volume-bound}
\sqrt{\frac{\det\Cov(\mu)}{\det\Cov(\nu)}}\ls C^n(1+D/n)^n.
\end{equation}
\end{lemma}

\begin{proof} If $D=0$, the two probability
measures are equal and the assertion is immediate. So we assume that $D>0$. We write $f_\nu=M e^{-V}$, where $M=\|f_\nu\|_\infty$ and $V\gr0$ is convex.
We use the upper-semicontinuous log-concave representative of the density, so that its maximum is attained at some $x_0$ and $V(x_0)=0$.
For $0<s<1$, convexity gives $V(x_0+s(x-x_0))\ls sV(x)$. Changing variables we get
\begin{equation}\label{eq:power-partition-bound}
Z_s:=\int_{\mathbb{R}^n} e^{-sV(x)}dx\ls s^{-n}\int_{\mathbb{R}^n} e^{-V(x)}dx=s^{-n}/M.
\end{equation}
If $\nu_s$ is the log-concave probability measure with density $Z_s^{-1}e^{-sV}$ then  $$0 \ls D_{\mathrm{KL}}(\mu\|\nu_s)=-h(\mu)+s\mathbb E_\mu V+\log Z_s.$$
All quantities are finite here because $h(\mu) \in \mathbb{R}$ by Theorem~\ref{th:log-concave-entropy-bounds}, while the finiteness of $D$ implies that $\mathbb E_\mu V<+\infty$. Therefore,
$$h(\mu)\ls s\mathbb E_\mu V+\log Z_s.$$
Set $H=h(\mu)+\log M$. Since $D=\mathbb E_\mu V-H$, using \eqref{eq:power-partition-bound} we obtain
$$(1-s)H\ls sD-n\log s.$$
Choose $s=(1+D/n)^{-1}$ and set $x=D/n$. Then
$$H\ls n+n(1+x^{-1})\log(1+x)\ls2n+n\log(1+x),$$
using $\log(1+x)\ls x$. Since $f_\nu\ls M$, we have $h(\nu)+\log M\gr0$, which proves \eqref{eq:entropy-difference-bound}.

Applying the lower bound of Theorem~\ref{th:log-concave-entropy-bounds} to $h(\mu)$ and the upper bound to $h(\nu)$ yields \eqref{eq:covariance-volume-bound}.
\end{proof}

For the remainder of this subsection we fix $u \in S^{n-1}$ and we work on the dual ray $\theta=ru$, $r>0$. Set $d=n-1$.

\begin{lemma}\label{lem:d-cov-ratio-1}
If $0<r_1<r_2$ we can write 
\begin{equation}\label{eq:determinant-orthogonal-decomposition}
\det\Sigma(r_iu)=a(r_i,u) \det S_i, \quad i=1,2,
\end{equation}
where $S_1$ and $S_2$ are $d \times d$ covariance matrices with 
$$\sqrt{\frac{\det S_1}{\det S_2}} \ls C^d \left(1+\frac{D_{\mathrm{KL}}(\mu_{r_1u}\|\mu_{r_2u})}{d}\right)^d.$$ 
\end{lemma}

\begin{proof}
For $i=1,2$ let $X_i$ have distribution $\mu_{r_iu}$ and set $a_i=a(r_i,u)$. After extending $u$ to an orthonormal basis $\{u,u_1,\ldots,u_d\}$ of $\mathbb{R}^n$, there exists an orthogonal transformation $U$ such that $$\widetilde{X}_i=UX_i=(Z_i,Y_i), \quad Z_i=\langle u,X_i\rangle.$$ Then $\widetilde{X}_i$ has covariance matrix $$H_i=\begin{pmatrix}a_i&c_i^T\\c_i&B_i\end{pmatrix}.$$ If we write $$S_i=B_i-\frac{c_ic_i^T}{a_i}$$ for the covariance matrix of $\,Y_i-(c_i/a_i)Z_i$, Schur's block complement formula gives $\det H_i=a_i\det S_i$. As determinants are invariant under orthogonal transformations, this is exactly \eqref{eq:determinant-orthogonal-decomposition}. We use the surjective linear map $$T(z,y)=y-\frac{c_2}{a_2}z, \quad z \in \mathbb{R}, \, y \in \mathbb{R}^d $$
and let $\nu_i$ be the distribution of $T\widetilde{X}_i$. Their covariance matrices satisfy
\begin{align*}
\Cov(\nu_2)&=S_2,\\ \Cov(\nu_1)&=S_1+a_1\left(\frac{c_1}{a_1}-\frac{c_2}{a_2}\right)\left(\frac{c_1}{a_1}-\frac{c_2}{a_2}\right)^T\succeq S_1.
\end{align*}
The measures $\nu_i$ are the push-forwards of $\mu_{r_iu}$ under the surjective linear map $T \circ U$ and thus are log-concave probability measures on $\mathbb R^d$. Using \eqref{eq:data-processing-inequality} we get $$D(\nu_1\Vert\nu_2)\ls D(\mu_{r_1u}\Vert\mu_{r_2u}).$$  
We apply Lemma~\ref{lem:entropy-volume-comparison} in dimension $d$. Since positive definite matrix order implies determinant order,
\begin{equation}\label{eq:schur-entropy-bound}
\sqrt{\frac{\det S_1}{\det S_2}}\ls\sqrt{\frac{\det\Cov(\nu_1)}{\det\Cov(\nu_2)}}
\ls C^d\left(1+\frac{D(\mu_{r_1u}\Vert\mu_{r_2u})}{d}\right)^d.
\end{equation}
This completes the proof.\end{proof}

In the notation of Lemma~\ref{lem:d-cov-ratio-1} we have the following:

\begin{lemma}\label{lem:d-cov-ratio-2}
If $0<r_1<r_2$ then $$\left(\frac{r_1}{r_2}\right)^d \sqrt{\frac{\det S_1}{\det S_2}} \ls C^d\left(1+\frac{\Delta}{d}\right)^d,$$ where $\Delta=F(r_2u)-F(r_1u)$ and $C>0$ is an absolute constant. 
\end{lemma}

\begin{proof}
If we set $$m(r)=\langle\nabla \Lambda(ru),u \rangle=\mathbb E_{\mu_{ru}}\langle u,X\rangle, \quad r>0$$ then $m'(r)=a(r,u)$. We also remind that $\frac{dF(ru)}{dr}=ra(r,u)$. Direct evaluation of relative entropy and integration by parts give
\begin{align*}
D(\mu_{r_1u}\Vert\mu_{r_2u}) &=\Lambda(r_2u)-\Lambda(r_1u)-(r_2-r_1)m(r_1)\\
 &=\int_{r_1}^{r_2}(r_2-s)a(s,u)\,ds \\
 &\ls\left(\frac{r_2}{r_1}-1\right)\int_{r_1}^{r_2}s a(s,u)\,ds= \left(\gamma^{-1}-1\right)\Delta, 
\end{align*}
where $\gamma=r_1/r_2$. Combining this with Lemma~\ref{lem:d-cov-ratio-1} we get 
\begin{align*}\label{eq:scaled-covariance-bound}
\left(\frac{r_1}{r_2}\right)^d\sqrt{\frac{\det S_1}{\det S_2}} &\ls C^d\left[\gamma+(1-\gamma)\frac\Delta d\right]^d\\
 &\ls C^d(1+\Delta/d)^d,
\end{align*}
because $0<\gamma<1$.
\end{proof}

\subsection{Comparison of ray densities and tails}

In this subsection we prove Theorem~\ref{th:near-linear}. We first show the corresponding result for the modified wet part $\widetilde w_\mu$.
\begin{proposition}\label{prop:modified-wet-part}
If $\lambda >1 $ and $0<t\ls e^{-1}/\lambda$ then $$\widetilde w_\mu(\lambda t) \ls C^n \lambda\left(1+\frac{\log \lambda}{n-1}\right)^{n-1} \widetilde w_\mu(t),$$ where $C>0$ is an absolute constant.
\end{proposition}

\begin{proof}
Recall Proposition~\ref{prop:dual-polar} and fix $u \in S^{n-1}$. By Lemma~\ref{lem:G-main-properties} we may define $F_u(g)$ as the inverse of $G(F,u)$. 
We make the change of variables $G(F,u)=g$, which gives $$dF=\frac{dg}{G'(F_u(g),u)}.$$
Observe that for every $T>0$ 
\begin{equation}\label{eq:tilde-wet-part-formula}\widetilde w_\mu(e^{-T})=\int_{S^{n-1}}M_u(T)\,d\sigma(u), \quad M_u(T)=\int_{T}^{\infty} k_u(g)\, dg, \quad k_u(g)=\frac{p_u(F_u(g))}{G'(F_u(g),u)}.
\end{equation}
Let $G_2 > G_1 \gr 1$. Then we can find $0<r_1<r_2<\ell_\Omega(u)$ such that if we set $F_i=F(r_iu)$ we have $G(F_i,u)=G_i$ for $i=1,2$. Let $h=G_2-G_1$ and $\Delta=F_2-F_1$. Using the decomposition of Lemma~\ref{lem:d-cov-ratio-1} in the formulas of \eqref{eq:Kn}, \eqref{eq:ray-density-entropic} and \eqref{eq:Gdef} we obtain $$e^{G_i}\rho(F_i)=r_i^{n-1}\sqrt{\det S_i}\left(1+\frac{1}{J(r_i,u)}\right), \quad i=1,2.$$ Since $G_1,G_2 \gr 1$, Lemma~\ref{lem:G-main-properties} gives $$\frac{1+1/J(r_1,u)}{1+1/J(r_2,u)} \ls 3.$$
Applying Lemma~\ref{lem:d-cov-ratio-2} we deduce that 
$$\frac{e^{G_1}\rho(F_1)}{e^{G_2}\rho(F_2)}\ls3C^{n-1}(1+\Delta/(n-1))^{n-1}.$$
Since $h\gr c_0\Delta$ by the lower bound of $G'$ in Lemma~\ref{lem:G-main-properties}, the last equation yields $$\rho(F_1)\ls C_1^n e^h\left(1+\frac{h}{n-1}\right)^{n-1}\rho(F_2).$$ Let now $g \gr 1$ and $h > 0$. Taking into account \eqref{eq:ray-density-comparison} and the bounds of $G'$ from Lemma~\ref{lem:G-main-properties}, we can conclude that $$k_u(g) \ls C_2^n e^{h}\left(1+\frac{h}{n-1}\right)^{n-1}k_u(g+h).$$ If $T-h \gr 1$, integrating over $g \gr T-h$ we get $$M_u(T-h)\ls  C_2^n e^{h}\left(1+\frac{h}{n-1}\right)^{n-1} M_u(T)$$ and thus by equation \eqref{eq:tilde-wet-part-formula} $$\widetilde w_\mu\big(e^{-(T-h)}\big) \ls C_2^ne^{h}\left(1+\frac{h}{n-1}\right)^{n-1} \widetilde w_\mu\big(e^{-T}\big).$$ Setting $T=-\log t$ and $h=\log \lambda$ we conclude the proof.
\end{proof}

\begin{proof}[Proof of Theorem~$\ref{th:near-linear}$]
Theorem~\ref{th:pointwise} gives $$c_1\widetilde q\ls q_\mu\ls C_1\widetilde q$$ for some absolute constants $ C_1 \gr c_1>0$. If $\lambda t/c_1>e^{-1}$, 
we use $w_\mu(t)\gr t$ to obtain
$$w_\mu(\lambda t) \ls 1 \ls \frac{e}{c_1}\lambda t \ls\frac{e}{c_1}\lambda w_\mu(t),$$ which is stronger than the required estimate. So we may assume that $\lambda t/c_1\ls e^{-1}$. Set $K=C_1/c_1$ and apply  Proposition~\ref{prop:modified-wet-part} at $T=t/C_1$ with multiplier 
$L=K\lambda$: 
\begin{align*}
w_\mu(\lambda t)&\ls \widetilde{w}_\mu(\lambda t/c_1)=\widetilde w_\mu(LT) \ls C^n L\left(1+\frac{\log L}{n-1}\right)^{n-1}\widetilde w_\mu(T)\\
&\ls C^nL\left(1+\frac{\log L}{n-1}\right)^{n-1} w_\mu(C_1T)\\
&\ls C_2^n \lambda \left(1+\frac{\log K+\log \lambda}{n-1}\right)^{n-1} w_\mu(t).
\end{align*} The inequality $$\left(1 +\frac{x+y}{n-1}\right)^{n-1}\ls e^x \left(1 +\frac{y}{n-1}\right)^{n-1}, \quad x,y \gr 0$$ concludes the proof of \eqref{eq:near-linear-log}. Let now $\varepsilon \in (0,1)$. Since $$\sup_{x \gr 0}(1+x)e^{-\varepsilon x}=\frac{e^{-(1-\varepsilon)}}{\varepsilon} \ls \frac{1}{\varepsilon},$$ setting $x=\log \lambda/(n-1)$ we get $$\left(1+\frac{\log \lambda}{n-1}\right)^{n-1} \ls \varepsilon^{-(n-1)}\lambda^\varepsilon.$$ Inserting this into \eqref{eq:near-linear-log} we complete the proof.
\end{proof}

\section{Counterexample for exponent one}\label{sec:endpoint-counterexample}

In this section we show that the factor multiplying $\lambda$ in Theorem~\ref{th:near-linear} cannot be replaced by a constant. Let $n \gr 2$. We work for $a \gr a_n$ for some sufficiently large constant $a_n>0$. To simplify the notation, throughout this section we omit the dependence of the functions and measures on $n$. The constants $A_n,B_n, C_n$ etc. depend only on $n$, and their values may change from line to line. 
Define $$\psi_a(p)=\begin{cases}
\max\{ |p|^2/2, \, na\|p\|_\infty-na^2/2\} &, \ p \in [-a,a]^n\\
+\infty &, \text{ otherwise}. 
\end{cases}
$$
Observe that $\psi_a$ is convex and set $V_a=\psi_a^\ast$, which is also convex.

\begin{lemma}\label{lem:V-basic-properties}
For every $x \in \mathbb{R}^n$ \begin{equation}\label{eq:Va-bounds}
\max\{a\|x\|_1-na^2/2, \, H_a(x)\}\ls V_a(x) \ls \frac{|x|^2}{2},\end{equation}
where $$H_a(x)=\begin{cases}
\frac{|x|^2}{2} &, \, |x| \ls a/2\\
\frac{a}{2}|x|-\frac{a^2}{8} &, \, |x| \gr a/2.
\end{cases}$$
In addition, we have that
\begin{equation}\label{eq:Va-equality} V_a(x)=a\|x\|_1-\frac{na^2}{2}, \quad \|x\|_1 \gr na.
\end{equation}
\end{lemma}

\begin{proof}
If $x \in \mathbb{R}^n$ we can see that $$V_a(x) \ls \sup_{p \in [-a,a]^n} \{\langle p,x \rangle-|p|^2/2\} \ls \sup_{p \in \mathbb{R}^n} \,\{\langle p,x \rangle-|p|^2/2\}=\frac{|x|^2}{2}. $$ On the other hand, choose $p_x\in\{-a,a\}^n$ by setting its $j$-th coordinate equal to $a$ if $x_j \gr 0$, and to $-a$ otherwise. Then $\psi_a(p_x)=na^2/2$ and $$V_a(x) \gr \langle p_x,x \rangle -\psi_a(p_x) =  a\|x\|_1-na^2/2.$$
For the other lower bound, if $|x| \ls a/2$ then $\|x\|_\infty \ls a/2$ and thus $$V_a(x) \gr \sup_{p \in [-a/2,a/2]^n} \{\langle p,x \rangle -\psi_a(p)\}=\sup_{p \in [-a/2,a/2]^n}\{\langle p,x \rangle -|p|^2/2\}=\frac{|x|^2}{2}.$$ If $|x| \gr a/2$ and $p_x=\frac{a}{2}\frac{x}{|x|}$ then $\|p_x\|_\infty \ls a/2$ and $$V_a(x) \gr \langle p_x , x
\rangle -\psi_a(p_x)=\frac{a}{2}|x|-\frac{a^2}{8}.$$ To conclude, if $\|x\|_1\gr na$, then 
\begin{align*}
V_a(x)&\ls\sup_{p\in[-a,a]^n}\{\langle p,x\rangle-na\|p\|_\infty+na^2/2\}\\
&\ls\sup_{0\ls m\ls a}\{m(\|x\|_1-na)+na^2/2\}
=a\|x\|_1-na^2/2.
\end{align*}
Combining this with the corresponding lower bound proves \eqref{eq:Va-equality}.
\end{proof}

Set now $$Z_a=\int_{\mathbb{R}^n }e^{-V_a(x)}\, dx, \quad f_a=\frac{1}{Z_a}e^{-V_a}$$ Then $f_a$ is a log-concave density on $\mathbb{R}^n$. We denote by $\mu_a$ the log-concave probability measure with density $f_a$. We note that $\mu_a$ is symmetric, since $\psi_a$ is even. Write also  $\phi(r)=(2\pi)^{-1/2}e^{-r^2/2}$ and $\overline\Phi(r)=\int_r^\infty\phi(s)\,ds$.

\begin{lemma}
If $a$ is large enough then \begin{equation}\label{eq:Za-bounds}
(2\pi)^{n/2} \ls Z_a \ls (2\pi)^{n/2}+C_n a^{n-2}e^{-a^2/8},
\end{equation} for some constant $C_n>0$. In particular, \begin{equation}\label{eq:Za-approx-constant}
Z_a \approx_n 1.
\end{equation}
Moreover, if $B_a=\frac{a}{2}B_2^n$ then \begin{equation}\label{eq:ball-mass-upper}
\mu_a(B_a^c) \ls \overline \Phi(a/4).
\end{equation}
\end{lemma}

\begin{proof}
The upper bound in \eqref{eq:Va-bounds} gives the inequality $Z_a \gr (2\pi)^{n/2}$, while the corresponding lower bound yields $$Z_a  \ls (2\pi)^{n/2}+e^{a^2/8}\int_{{B_a^c}}e^{-a|x|/2} \,dx.$$ Using polar coordinates and making the change of variables $t=ar/2$, we obtain $$\int_{{B_a^c}}e^{-a|x|/2}\, dx=n\omega_n\int_{a/2}^{\infty} r^{n-1}e^{-ar/2}\, dr= n\omega_n\left(\frac{2}{a}\right)^n\int_{a^2/4}^{\infty} t^{n-1}e^{-t}\, dt.$$ The fact that $\int_{x}^{\infty}t^{n-1}e^{-t}\,dt\sim x^{n-1}e^{-x}$, as $x \to \infty$, implies $$Z_a \ls (2\pi)^{n/2}+C_n a^{n-2}e^{-a^2/8},$$ for $a$ large enough. The last calculation combined with the lower bound for $Z_a$ shows that $$\mu_a(B_a^c) \ls A_n a^{n-2}e^{-a^2/8},$$ while integration by parts gives \begin{equation}\label{eq:gaussian-tail-estimates}\frac{r}{r^2+1} \phi(r) \ls\overline \Phi(r) \ls\frac{\phi(r)}{r}, \quad r>0.
\end{equation}
Therefore, if $a$ is large enough, we have $$\frac{\mu_a(B_a^c)}{\overline\Phi(a/4)}\ls D_n a^{n-1} e^{-3a^2/32}<1.$$ 
\end{proof}

\begin{lemma}\label{lem:sa-definition-ration}
Set $r_a=a/4$ and $s_a=2\overline\Phi(r_a)$. Then we have that $$\frac{w_{\mu_a}(s_a)}{s_a}\gr c_na^{n-1}.$$
\end{lemma}

\begin{proof}
We will first show that \begin{equation}\label{eq:wet-part-lower}\{x: |x| \gr r_a\} \subseteq \{x:q_{\mu_a}(x) \ls s_a\}.\end{equation}
Indeed, if $|x| \gr r_a$ we set $u=x/|x|$ and consider the half-space $H_x=\{z: \langle z,u \rangle \gr r_a \}$ that contains $x$. Then, $$q_{\mu_a}(x) \ls \mu_a(H_x).$$ Using \eqref{eq:ball-mass-upper} we obtain $$\mu_a(H_x) \ls \mu_a(H_x\cap B_a)+\mu_a(B_a^c) \ls \mu_a(H_x\cap B_a) +\overline \Phi(r_a),$$ while the lower bounds of \eqref{eq:Va-bounds} and \eqref{eq:Za-bounds} give $$\mu_a(H_x \cap B_a)  \ls \gamma_n(H_x) =\overline\Phi(r_a),$$ where $\gamma_n$ is the standard Gaussian measure in $\mathbb{R}^n$. The last equality follows by rotational invariance of $\gamma_n$. Hence, $$q_{\mu_a}(x) \ls s_a.$$ Taking measures in \eqref{eq:wet-part-lower} yields $$\mu_a(\{x: |x| \gr r_a\}) \ls w_{\mu_a}(s_a).$$ Using the upper bound of \eqref{eq:Va-bounds} and \eqref{eq:Za-approx-constant}, as well as polar coordinates we get \begin{align*}    
\mu_a(\{x: |x| \gr r_a\}) & \gr  A_n \int_{r_a}^{r_a+\frac{1}{r_a}} r^{n-1}e^{-r^2/2}\, dr \gr A_n r_a^{n-2}e^{-(r_a+1/r_a)^2/2} \\
&\gr B_n r_a^{n-2}e^{-r_a^2/2}=C_na^{n-2}e^{-a^2/32},
\end{align*} if $a$ is large enough. The upper bound in \eqref{eq:gaussian-tail-estimates} completes the proof.
\end{proof}

Let  $\nu_a$  be the log-concave probability measure with density $$g_a(x)=\bigg(\frac{a}{2}\bigg)^ne^{-a\|x\|_1}.$$ Observe that $\nu_a$ is a product measure. The following lemma is a direct consequence of \eqref{eq:Va-bounds} and \eqref{eq:Va-equality}.

\begin{lemma}\label{lem:laplace-density-comparison}
There exists a constant $C_{a,n}>0$ such that $f_a \ls C_{a,n}g_a$ and $$f_a(x)=C_{a,n}g_a(x), \quad \|x\|_1 \gr na.$$ In particular, $$C_{a,n}=\left(\frac{2}{a}\right)^n\frac{e^{a^2n/2}}{Z_a}.$$
\end{lemma}

\begin{lemma}\label{lem:laplace-depth-comparison}
If $a$ is large enough and $\|x\|_\infty \gr 2na$ then $$q_{\mu_a}(x)=C_{a,n}q_{\nu_a}(x).$$
\end{lemma}

\begin{proof}
After normalization, for any probability measure $\mu$ we may write $$q_{\mu}(x)=\inf_{\|u\|_\infty=1} \mu(\{z: \langle z,u \rangle \gr \langle x,u \rangle\}).$$ We first claim that it suffices to restrict ourselves to half-spaces with $\langle u, x\rangle \gr na$. If $\|u\|_\infty=1$, we will obtain an  appropriate lower bound for the measure of $H=\{z:\langle z, u\rangle \gr na\}$, with respect to $\mu_a$ and $\nu_a$. Choose a coordinate $j$ with $|u_j|=1$. Consider the rectangle $R$ with $$\begin{cases}
-\frac{1}{2n} \ls z_i \ls \frac{1}{2n} &, i \neq j\\ 
na+1 \ls {\rm sgn(u_j)} z_j \ls na+2 &, i=j.
\end{cases} $$ Every point $z \in R$ satisfies $$\|z\|_1 \gr\langle u,z \rangle \gr na+1-\frac{n-1}{2n}>na.$$ Hence, $R$ lies in the half-space $\{z:\langle u,z \rangle \gr na \}$ and $$\nu_a(R)=\frac{1}{2}(1-e^{-a})(1-e^{-a/(2n)})^{n-1}e^{-na^2-a}.$$ Using Lemma~\ref{lem:laplace-density-comparison}, if $a$ is large enough, we deduce that \begin{equation}\label{eq:half-spaces-lower-bound}
\nu_a(H) \gr \frac{1}{4}e^{-na^2-a}, \quad \mu_a(H)\gr \frac{C_{a,n}}{4}e^{-na^2-a}.
\end{equation} If $|x_j|=\|x\|_\infty$ then the coordinate half-space $H_j=
\{z: {\rm sgn}(x_j)z_j \gr \|x\|_\infty \}$ passes through $x$ and has $\nu_a(H_j)=\frac{1}{2}e^{-a\|x\|_\infty}$ and $\mu_a(H_j)=\frac{C_{a,n}}{2}e^{-a\|x\|_\infty}$. Therefore, $$q_{\nu_a}(x) \ls \frac{1}{2}e^{-2na^2}, \qquad q_{\mu_a}(x) \ls \frac{C_{a,n}}{2}e^{-2na^2}.$$ For $a$ large enough, those upper bounds are smaller than the lower bounds in \eqref{eq:half-spaces-lower-bound}, which implies the claim. The proof now follows from Lemma~\ref{lem:laplace-density-comparison}, since every such half-space is contained in $\{z: \|z\|_1 \gr na\}$. 
\end{proof}

\begin{lemma}\label{lem:depth-lower-bounds-cube}
There exists a constant $A_n>0$ such that $$\inf_{\|x\|_\infty \ls 2na} q_{\mu_a}(x) \gr e^{-A_na^2}, \quad \inf_{\|x\|_\infty \ls 2na} q_{\nu_a}(x) \gr e^{-A_na^2},$$ if $a$ is large enough.
\end{lemma}

\begin{proof}
If $H$ is a half-space that passes through $x$, it contains half of the unit ball, centered at $x$. If $y$ is a point of this ball, then $$|y| \ls 2n^{3/2}a+1 \ls 3n^{3/2}a.$$ Therefore, the upper bound of \eqref{eq:Va-bounds} and \eqref{eq:Za-approx-constant} give $$f_a(y) \gr c_ne^{-C_na^2}.$$ The inequality $\|y\|_1 \ls \sqrt{n}|y|$ gives the lower bound $$g_a(y) \gr d_n e^{-B_na^2}.$$ Integrating the densities over this half-ball and taking infimum over $u$ concludes the proof.
\end{proof}

Set now $D_n=2\max\{A_n,1\}$, where $A_n$ is the constant from Lemma~\ref{lem:depth-lower-bounds-cube} and define \begin{equation}\label{eq:tau-a-definition}
t_a=e^{-D_na^2}, \quad \tau_a=\frac{t_a}{C_{a,n}}.
\end{equation}

\begin{lemma}\label{lem:wet-parts-ratio-comparison}
If $a$ is large enough then $$\frac{w_{\mu_a}(t_a)}{t_a} \ls \frac{w_{\nu_a}(\tau_a)}{\tau_a}.$$
\end{lemma}

\begin{proof}
If $a$ is large enough, then Lemma~\ref{lem:depth-lower-bounds-cube} and the definition of $D_n$ force $$W_{\mu_{a}}(t_a) \subseteq\{x: \|x\|_\infty \gr 2na\}. $$ Therefore, Lemma~\ref{lem:laplace-density-comparison} and Lemma~\ref{lem:laplace-depth-comparison} give $$w_{\mu_a}(t_a)=\int_{\{x:q_{\mu_a}(x) \ls t_a\}} f_a(x) \, dx \ls C_{a,n}\int_{\{x:q_{\nu_a}(x) \ls \tau_a\}} g_a(x) \, dx=C_{a,n}w_{\nu_{a}}(\tau_a). $$ Dividing by $t_a$ yields the proof.
\end{proof}
We will now estimate $w_{\nu_a}(t)$ for sufficiently small $t>0$.

\begin{lemma}\label{lem:laplace-wet-part-legendre-bound}
There exists a constant $t_0 \in (0,1)$ such that for every $a>0$ and every $0<t<t_0$ we have that $$w_{\nu_a}(t) \ls \nu_1(\{x:\Lambda_{\nu_1}^\ast(x) \gr b(t) \}),$$ where $b(t)$ is a function satisfying \begin{equation}\label{eq:b-asymptotics}\frac{1}{2}\log \frac{1}{t}\ls b(t) \ls \log \frac{1}{t}, \quad e^{-b(t) }\ls 2et\log \frac{1}{t}\end{equation}
for every $0<t<t_0$.
\end{lemma}

\begin{proof}
By affine invariance we can see that $w_{\nu_a}(t)=w_{\nu_1}(t)$ for every $a,t>0$. Let $0<t<e^{-1}$. Using \eqref{eq:BC-inequality} for $x \in W_{\nu_1}(t)$ and $\varepsilon=1/\log(1/t)$ we get $$\Lambda_{\nu_1}^\ast(x) \gr b(t) ,$$ where $$b(t)=\log(1/t)-\log\log(1/t)-1-\log2+\log2/\log(1/t).$$ For $0<t<e^{-1}$ we have $\log(1/t)>1$, which implies that $$\log(1/t)-\log\log(1/t)-1-\log 2<b(t) <\log(1/t).$$ 
Composing the left-hand side inequality with the exponential function we get the right-hand side inequality of \eqref{eq:b-asymptotics}. The lower bound for $b(t)$ follows from the fact that $b(t)/\log(1/t) \to 1$, as $t \to 0^+$. 
\end{proof}

\begin{lemma}\label{lem:laplace-legendre-tail-bound}
There exist constants $T_n,B_n>0$ such that for every $t \gr T_n$ we have $$\nu_1(\{x: \Lambda_{\nu_1}^\ast(x) \gr t \}) \ls B_n e^{-t}\frac{(\log t)^{n-1}}{t}.$$
\end{lemma}

\begin{proof}
Let $\nu$ be the one-dimensional marginal of $\nu_1$. Since $\nu_1$ is a product measure we have that $$\Lambda_{\nu_1}^\ast(x)=\sum_{i=1}^{n} \Lambda_\nu^\ast(x_i).$$ A direct calculation shows that $$\Lambda_\nu^\ast(x)=\sqrt{1+x^2}-1+\log\biggl(\frac{2}{1+\sqrt{1+x^2}}\biggr), \quad (\Lambda^\ast_\nu)^\prime(x)=\frac{x}{1+\sqrt{1+x^2}},$$ for every $x \in \mathbb{R}$. Observe that $\Lambda_\nu^\ast$ is even and has an inverse on $(0,\infty)$. If $X$ is a random variable distributed according to $\nu$, then the non-negative random variable $Y=\Lambda^\ast_\nu(X)$ has density $$f_Y(y)=\frac{e^{-x}}{(\Lambda^\ast_\nu)^\prime(x)}=e^{-x}\frac{1+\sqrt{1+x^2}}{x}, \quad x=(\Lambda_\nu^\ast)^{-1}(y), \quad y>0.$$ Since $(\Lambda^\ast_\nu)^\prime(x) \ls \frac{x}{1+x}$, integrating from $0$ to $x$ gives $$\Lambda^\ast_{\nu}(x) \ls x-\log(1+x).$$ Thus \begin{equation}\label{eq:symm-exp-first-bound}
e^{-x}(1+x) \ls e^{-\Lambda_\nu^\ast(x)}.
\end{equation}
Using \eqref{eq:symm-exp-first-bound} and the inequality $$1+\sqrt{1+x^2} \ls 2(1+x)$$ we get $$f_Y(y) \ls 2e^{-x}\frac{1+x}{x}\ls \frac{2e^{-\Lambda_\nu^\ast(x)}}{x}=\frac{2e^{-y}}{x}.$$ We will now show that $$\Lambda^\ast_\nu(x)(1+\Lambda^\ast_\nu(x)) \ls 2x^2$$ for every $x \gr 0$, which implies that \begin{equation}\label{eq:density-upper-bound}
f_Y(y) \ls 2\sqrt{2}\frac{e^{-y}}{\sqrt{y(y+1)}}, \quad y>0.
\end{equation}
Indeed, if $0 \ls x \ls 1$ then $(\Lambda_\nu^\ast)^\prime(x) \ls \frac{x}{2}$ and integration from $0$ to $x$ gives $\Lambda^\ast_\nu(x) \ls \frac{x^2}{4}$. The desired inequality follows. For $x \gr 1$ we use the trivial inequality $(\Lambda^\ast_\nu)^\prime(x) \ls 1$, which implies that $\Lambda_\nu^\ast(x) \ls x$. Therefore, since $x \gr 1$ $$\Lambda^\ast_\nu(x)(\Lambda^\ast_\nu(x)+1) \ls x(x+1) \ls 2x^2.$$ If $X_1,\ldots,X_n$ are i.i.d random variables distributed according to $\nu$, the convolution formula implies that the density of $S_n=\sum_{i=1}^n\Lambda^\ast_\nu(X_i)$ satisfies \begin{equation}\label{eq:sum-density-upper-bound}f_{S_n}(s) \ls C^n e^{-s} h^{\ast n}(s) , \quad s>0,
\end{equation}
where $C=2\sqrt{2}$ and $$h(y)=\frac{1}{\sqrt{y(y+1)}}, \quad y>0.$$ We will show by induction that \begin{equation}\label{eq:convolution-upper-bound}h^{\ast n}(s) \ls C_n \frac{(\log s)^{n-1}}{s}, \quad s \gr s_n,
\end{equation}for some $s_n \gr 1$. To begin with, it is obvious that for $s >0$ we have $h(s) \ls 1/s$. Assume that for some $k$ we have $$h^{\ast k}(s) \ls C_k \frac{(\log s)^{k-1}}{s}, \quad s \gr s_k,$$ where $s_k \gr 1$. Let $s \gr 2s_{k}$ and write $$h^{\ast (k+1)}(s)=\int_{0}^{s/2} h^{\ast k}(u)h(s-u)\, du+ \int_{s/2}^{s} h^{\ast k}(u)h(s-u)\, du.$$ For the first integral $s-u \gr s/2$ implies $h(s-u) \ls 2/s$. By the induction hypothesis,  $$\int_{0}^{s/2} h^{\ast k}(u)h(s-u)\, du \ls \frac{2}{s}\bigg[D_k+C_k\int_{s_k}^{s/2} \frac{(\log u)^{k-1}}{u} \, du \bigg] \ls A_k \frac{(\log s)^k}{s},$$ where we used the fact that $h^{\ast k}$ is integrable on $[0,s_k]$, since $h$ is integrable and non-negative on the same interval. For the second integral, the induction hypothesis followed by a change of variables, gives $$\int_{s/2}^{s} h^{\ast k}(u)h(s-u)\, du \ls \frac{2C_k(\log s)^{k-1}}{s}\int_{0}^{s/2} h(u)\, du \ls B_k\frac{(\log s)^k}{s},$$ because $$\int_{0}^{s/2} h(u) \, du \ls \int_{0}^{1} \frac{1}{\sqrt{u}}\, du +\int_{1}^{s/2} \frac{1}{u} \, du \ls C\log s.$$ These two estimates give the desired upper bound.
    
Since the function $$s \longrightarrow \frac{(\log s)^{n-1}}{s}$$ is decreasing for $s \gr e^{n-1}$, integration for $t \gr \max\{e^{n-1},s_n\}$ and the inequalities \eqref{eq:sum-density-upper-bound} and \eqref{eq:convolution-upper-bound} conclude the proof.
\end{proof}

\begin{lemma}\label{lem:ta-ratio-bound}
There exist constants $t_n,C_n>0$ such that for every $0<t<t_n$ and every $a>0$  \begin{equation}\label{eq:laplace-wet-part-ratio}
\frac{w_{\nu_a}(t)}{t}\ls C_n \bigg(\log \log\frac{1}{t}\bigg)^{n-1}.\end{equation}
In particular, if $a$ is large enough then $$\frac{w_{\mu_a}(t_a)}{t_a}\ls K_n (\log a)^{n-1},$$ for some constant $K_n>0$.
\end{lemma}

\begin{proof}
By Lemma~\ref{lem:laplace-wet-part-legendre-bound} we can find $t_n \in (0,t_0)$ such that for every $0<t<t_n$ we have $b(t)>T_n$. Therefore, Lemma~\ref{lem:laplace-wet-part-legendre-bound} and Lemma~\ref{lem:laplace-legendre-tail-bound} yield  $$\frac{w_{\nu_a}(t)}{t}\ls B_n \frac{e^{-b(t)}}{tb(t)}(\log b(t))^{n-1}.$$ The inequalities in \eqref{eq:b-asymptotics} give $$\frac{e^{-b(t)}}{tb(t)} (\log b(t))^{n-1} \ls 4e \left(\log\log\frac{1}{t}\right)^{n-1}.$$ This proves \eqref{eq:laplace-wet-part-ratio}. Taking into account \eqref{eq:tau-a-definition} and the definition of $C_{a,n}$ in Lemma~\ref{lem:laplace-density-comparison}, we see that if $a$ is large enough then $\tau_a<t_n$ and $$\tau_a \gr e^{-L_na^2}$$ for some constant $L_n>0$. Then Lemma~\ref{lem:wet-parts-ratio-comparison} and \eqref{eq:laplace-wet-part-ratio} conclude the proof.
\end{proof}

\begin{proof}[Proof of Theorem~$\ref{th:general-failure}$]
For $a$ sufficiently large, the symmetric log-concave measure $\mu_a$ constructed above satisfies $0<t_a<s_a<1$, by the definitions in Lemma~\ref{lem:sa-definition-ration} and \eqref{eq:tau-a-definition}, together with \eqref{eq:gaussian-tail-estimates}. Set $\lambda_a=s_a/t_a$. Then $\lambda_a>1$ and $0<t_a<1/\lambda_a$. Lemmas~\ref{lem:sa-definition-ration} and~\ref{lem:ta-ratio-bound} give
$$\frac{w_{\mu_a}(\lambda_a t_a)}{\lambda_a w_{\mu_a}(t_a)}=\frac{w_{\mu_a}(s_a)/s_a}{w_{\mu_a}(t_a)/t_a}
\gr \frac{c_n}{K_n}\left(\frac{a}{\log a}\right)^{n-1}\longrightarrow\infty\qquad(a\to\infty).$$
Given any $C_n>0$, choosing $a$ sufficiently large proves the theorem.
\end{proof}

\bigskip 

\noindent\textbf{AI Usage Statement.} The authors acknowledge the assistance of ChatGPT (OpenAI) at the final stage of the preparation of this work: this included proof checking, corrections and suggestions on the presentation.

\bigskip 

\medskip 


\bigskip 

\medskip

\noindent\textbf{2020 Mathematics Subject Classification.} Primary 60D05, 52A23; Secondary 52A22, 60E15, 52A40, 94A17.

\smallskip 

\noindent\textbf{Keywords.} Log-concave probability measures, Tukey half-space depth, wet-part doubling, Cram\'{e}r transform, exponential tilting, relative entropy, covariance determinants, random polytopes.

\bigskip

\bigskip

\bigskip

\noindent \textsc{Apostolos \ Giannopoulos}: School of Applied Mathematical and Physical Sciences, National Technical University of Athens, Department of Mathematics, Zografou Campus, GR-157 80, Athens, Greece.

\smallskip

\noindent \textit{E-mail:} \texttt{apgiannop@math.ntua.gr}

\bigskip

\noindent \textsc{Minas \ Pafis}: Department of Mathematics, National and Kapodistrian University of Athens, Panepistimioupolis 157-84,
Athens, Greece.

\smallskip

\noindent \textit{E-mail:} \texttt{mipafis@math.uoa.gr}

\end{document}